\title{ ~~\\ On the density of primes in arithmetic progression
having a prescribed primitive root}
\author{Pieter Moree}
\documentclass[12pt]{article}
\usepackage{amssymb, latexsym, amsfonts}
\def\@ptsize{2}
\newtheorem{Thm}{Theorem}
\newtheorem{lem}{Lemma}

\newtheorem{Prop}{Proposition}
\newcommand{\qed}{\hfill $\Box$}

\begin{document}
\date{}
\maketitle
{\def\thefootnote{}
\footnote{{\it 1991 Mathematics Subject Classificiation}. 11R45, 11A07,
11N69.}}
\begin{abstract}
\noindent
Let $g\in \mathbb Q$ be not $-1$ or a square.
Let ${\cal P}_g$ denote the set of primes $p$ such that $g$ is a primitive root
mod $p$. Let
$1\le a\le f,~(a,f)=1$. Under the Generalized Riemann Hypothesis (GRH)
it can be shown that the
set of primes $p\in {\cal P}_g$
with $p\equiv a({\rm mod~}f)$ has a natural density.
In this note this density is explicitly evaluated. This generalizes a
classical result of Hooley.
\end{abstract}
\section{Introduction}
Let $G$ be the set of non-zero rational
numbers $g$ such that $g\ne -1$ and $g$ is
not a square of a rational number. For
arbitrary $g\in \mathbb Q^*$ let ${\cal P}_g$ denote the
set of primes $p$ such that $g$ is a primitive root mod $p$.
Clearly a necessary condition for ${\cal P}_g$ to be infinite is
that $g\in G$. That this is also a sufficient condition was
conjectured in 1927 by Emil Artin.
There is no $g\in G$ for
which this has been proved, however, Heath-Brown
\cite{HB} in a classical paper
established a result which implies, for example, that there are at most
two primes $q$ for which ${\cal P}_q$ is finite. In 1967 Hooley
\cite{Hooley}
established Artin's primitive root conjecture under the assumption of
the Generalized Riemann Hypothesis (GRH). Moreover, he showed
that under that assumption the set ${\cal P}_g$ has a natural density,
which he evaluated (his result is
Theorem \ref{main} below with $a=f=1$ and
$g\in G\cap \mathbb Z$).
It turns out that this density is equal to a positive rational
number times the Artin constant $A$,
with
$$A=\prod_p\left(1-{1\over p(p-1)}\right)
\approx .373955813619202288054728.$$
(Throughout this note the notation $p$ is used to indicate rational primes.)\\
\indent In connection with his study of Euclidian number fields, Lenstra
\cite{L}
considered the distribution over
arithmetic progressions of the primes in ${\cal P}_g$.
Let ${\cal P}_{a,f,g}$ denote the set of primes $p$ such that $g$ is a primitive
root mod $p$ and $p\equiv a({\rm mod~}f)$. From Lenstra's work it follows
that, under GRH,  ${\cal P}_{a,f,g}$ has a natural density.
\begin{Thm} {\rm \cite{L}}.
\label{lens}
Put $\zeta_m=e^{2\pi i/m}$. Let $1\le a\le f,~(a,f)=1$. Let
$\sigma_a$ be the automorphism of $\mathbb Q(\zeta_f)$ determined by
$\sigma_a(\zeta_f)=\zeta_f^a$. Let $c_a(n)$ be $1$ if the restriction
of $\sigma_a$ to the field $\mathbb Q(\zeta_f)\cap \mathbb Q(\zeta_n,g^{1/n})$
is the identity and $c_a(n)=0$ otherwise. Put
$$\delta(a,f,g)=\sum_{n=1}^{\infty}{\mu(n)c_a(n)\over
[\mathbb Q(\zeta_f,\zeta_n,g^{1/n}):\mathbb Q]}.$$
Then, assuming GRH, we have
$$\pi_{g}(x;f,a)=\delta(a,f,g){x\over \log x}+
O\left({x\log \log x\over \log^2x}\right),$$
where $\pi_{g}(x;f,a)$ denotes the number of primes $p\le x$ that
are in ${\cal P}_{a,f,g}$.
\end{Thm}
In view of the apparent arithmetical complexity of Lenstra's formula, the
following
relatively simple expression for $\delta(a,f,g)$ may come as a
bit of a surprise.
\begin{Thm}
\label{main}
Let $g\in G$ and let
$h$ be the largest integer such that $g$ is an $h$-th power.
Let $\Delta$ denote the discriminant of the quadratic
field $\mathbb Q(\sqrt{g})$. Let $1\le a\le f,~(a,f)=1$.
Let $b=\Delta/(f,\Delta)$. If $b$ is odd, put
$\gamma=(-1)^{(b-1)/2}(f,\Delta)$. Put
$$A(a,f,h)=\prod_{p|(a-1,f)}(1-{1\over p})
\prod_{p\nmid f\atop p|h}(1-{1\over p-1})\prod_{p\nmid f\atop p\nmid h}
\left(1-{1\over p(p-1)}\right)$$
if $(a-1,f,h)=1$ and $A(a,f,h)=0$ otherwise. Then
$$\delta(a,f,g)={A(a,f,h)\over \varphi(f)}\left(1+
({\gamma\over a}){\mu(2|b|)\over
\prod_{p|b,~p|h}(p-2)\prod_{p|b,~p\nmid h}(p^2-p-1)}\right).$$
Here $(\gamma/a)$ denotes the Kronecker symbol.
\end{Thm}
In case $b$ is even, the Kronecker symbol $(\gamma/a)$ is not defined,
then however $\mu(2|b|)=0$ and the term involving $\mu(2|b|)$ is taken to
be zero. Note that $g$
can be uniquely written as $g=g_1g_2^2$, with $g_1$ a squarefree
integer
and $g_2\in \mathbb Q$. Then $\Delta=g_1$
if $g_1\equiv 1({\rm mod~}4)$ and $\Delta=4g_1$ otherwise. We see that
$b$ is odd if and only if $g_1\equiv 1({\rm mod~}4)$
or $g_1\equiv 2({\rm mod~}4)$ and $8|f$
or $g_1\equiv 3({\rm mod~}4)$ and $4|f$. On using the
properties of the Kronecker symbol mentioned in Section
\ref{algfacts} and quadratic reciprocity for the Jacobi symbol,
Theorem \ref{main} can be formulated
in terms of $g_1$. (The odd part of $0\ne m\in \mathbb Z$
is $m/2^e$ with $2^e|m$ and $2^{e+1}\nmid m$.)
\begin{Thm}
\label{main1}
Let $a,f,g,h$ and $A(a,f,h)$ be
as in Theorem {\rm \ref{main}}. Let $\tilde g_1$ and $\tilde a$ denote the
odd parts of $g_1$ and $a$, respectively. Let $\beta=g_1/(g_1,f)$. We have
$$\delta(a,f,g)={A(a,f,h)\over \varphi(f)}\left(1-
({a\over (f,g_1)}){(-1)^{{\tilde g_1-1\over 2}
{\tilde a-1\over 2}}\mu(|\beta|)\over
\prod_{p|\beta,~p|h}(p-2)\prod_{p|\beta,~p\nmid h}(p^2-p-1)}\right)$$
in case $g_1\equiv 1({\rm mod~}4)$ or $g_1\equiv 2({\rm mod~}4)$ and
$8|f$ or $g_1\equiv 3({\rm mod~}4)$ and $4|f$ and
$$\delta(a,f,g)={A(a,f,h)\over \varphi(f)},$$
otherwise. Here $(./.)$ denotes the Jacobi symbol with the stipulation
that $(a/2)=(-1)^{(a^2-1)/8}$.
\end{Thm}
From Theorem \ref{main} various known results can be rather easily deduced.
That will be the subject of
Section \ref{applic}. The remaining sections
are devoted to proving Theorem \ref{main}.\\
\indent The author thanks F. Lemmermeyer and P. Stevenhagen for some
helpful remarks. The research presented here was mostly carried out
in the pleasant and inspiring atmosphere of the Max-Planck-Institut
in Bonn.

\section{Some facts from algebraic number theory}
\label{algfacts}
Although Theorem \ref{lens} suggests differently, the problem
of evaluating $\delta(a,f,g)$ really only
involves cyclotomic fields, quadratic fields
and their composita. In this section some facts
concerning these fields relevant for the proof of
Theorem \ref{main} are discussed. We start by recalling some properties
of the Kronecker symbol, a rarely covered topic in books on number theory
(in contrast to the Jacobi symbol).\\
\indent The Kronecker symbol $(c/d)$ is defined for $c\in \mathbb Z,~
c\equiv 0({\rm mod~}4)$ or $c\equiv 1({\rm mod~}4)$, $c$ not a square, and
$d\ge 1$ an integer. If $c$ and $d$ are such that
both the Kronecker and the Jacobi symbols are defined, then these symbols
coincide. If $(c,d)>1$, then $(c/d)=0$. If $(c,d)=1$, then $(c/d)=\pm 1$.
If $d_1,d_2\ge 1$, then $(c/d_1d_2)=(c/d_1)(c/d_2)$. If $c$ is odd, then
$(c/2)=$ Jacobi symbol $(2/|c|)$. If $d_1\equiv d_2({\rm mod~}|c|)$,
then $(c/d_1)=(c/d_2).$\\
\indent The following
lemma allows one to determine all quadratic subfields of a given cyclotomic
field (for a proof see e.g. \cite{W}, p. 163).
\begin{lem}
\label{classical}
Let $\mathbb Q(\sqrt{d})$ be a quadratic
field of discriminant $\Delta_d$.
Then
the smallest cyclotomic field containing $\mathbb Q(\sqrt{d})$ is
$\mathbb Q(\zeta_{|\Delta_d|})$.
\end{lem}
Consider the cyclotomic field $\mathbb Q(\zeta_f)$. There are $\varphi(f)$
distinct automorphisms given by $\sigma_a(\zeta_f)=\zeta_f^a$, with
$1\le a\le f$ and $(a,f)=1$. We need
to know when the restriction of such an automorphism to a given
quadratic subfield of $\mathbb Q(\zeta_f)$ is the identity. In this
direction we have
\begin{lem}
\label{twee}
Let $\mathbb Q(\sqrt{d})\subseteq \mathbb Q(\zeta_f)$ be a quadratic
field of discriminant $\Delta_d$.
We have $\sigma_a|_{\mathbb Q(\sqrt{d})}=$id if and only if
$(\Delta_d/a)=1$, where $(./.)$ denotes the Kronecker symbol.
\end{lem}
{\it Proof}. We have
$$\sigma_a(\sqrt{d})=\left({\mathbb Q(\zeta_f)/\mathbb Q\over a}\right)\sqrt{d}
=\left({\mathbb Q(\sqrt{d})/\mathbb Q\over a}\right)\sqrt{d}=
({\Delta_d\over a})\sqrt{d},$$
where the latter reciprocity symbol is the Kronecker symbol and the other two are
Artin symbols. \qed

\noindent {\bf Remark}. It is also possible to prove
Lemma \ref{twee} using quadratic reciprocity and properties
of Gauss sums.\\
\indent The next result can be proved by a trivial generalization of an
argument given by Hooley \cite{Hooley}, pp. 213-214.
\begin{lem}
\label{cycdegree}
Let $g\in G$ and let $h$ be the
largest positive integer such that $g$ is
an $h$-th power. Let $\Delta$ denote the discriminant of
$\mathbb Q(\sqrt{g})$.
Let $k$ and $r$ be natural numbers such that
$k|r$ and $k$ is squarefree. Put $k_1=k/(k,h)$ and
$n(k,r)=[\mathbb Q(\zeta_r,g^{1/k}):\mathbb Q].$ Then
\begin{itemize}
\item[i)] if $k$ is odd, $n(k,r)=k_1\varphi(r);$
\item[ii)] if $k$ is even and $\Delta\nmid r,$ $n(k,r)=k_1\varphi(r);$
\item[iii)] if $k$ is even and $\Delta|r,$ $n(k,r)=k_1\varphi(r)/2.$
\end{itemize}
\end{lem}
Notice that $h$ is odd, since by assumption $g$ is not a square.
The next lemma together with Lemma \ref{twee} allows one to compute $c_a(n)$.
\begin{lem}
\label{gemier}
Let $g\in G$ and let $\Delta$ denote
the discriminant of $\mathbb Q(\sqrt{g})$. Let $n\ge 1$ be squarefree and
$f\ge 1$ be arbitrary. Put $b=\Delta/(f,\Delta)$. If $b$ is odd, then
$$\mathbb Q(\zeta_f) \cap \mathbb Q(\zeta_n,g^{1/n})=\mathbb Q(\zeta_{(f,n)},
\sqrt{(-1)^{(b-1)/2}(f,\Delta)})$$
if $n$ is even, $\Delta\nmid n$ and $\Delta|{\rm lcm}(f,n)$ and
$$\mathbb Q(\zeta_f) \cap \mathbb Q(\zeta_n,g^{1/n})=\mathbb Q(\zeta_{(f,n)})$$
otherwise. If $b$ is even, then $\mathbb Q(\zeta_f) \cap
\mathbb Q(\zeta_n,g^{1/n})=\mathbb Q(\zeta_{(f,n)}).$
\end{lem}
{\it Proof}. On noting that  $\varphi((f,n))\varphi({\rm lcm}(f,n))=
\varphi(f)\varphi(n)$ and that
$$[\mathbb Q(\zeta_f,\zeta_n,g^{1/n}):\mathbb Q]={[\mathbb Q(\zeta_f):\mathbb Q(\zeta_{(f,n)})][\mathbb Q(\zeta_n,g^{1/n}):\mathbb Q]\over
[\mathbb Q(\zeta_f)\cap \mathbb Q(\zeta_n,g^{1/n}):\mathbb Q(\zeta_{(f,n)})]}$$
it easily follows, using Lemma \ref{cycdegree},
that
\begin{equation}
\label{graad}
[\mathbb Q(\zeta_f)\cap \mathbb Q(\zeta_n,g^{1/n}):\mathbb Q(\zeta_{(f,n)})]=2
\end{equation}
if $n$ is even, $\Delta\nmid n$ and $\Delta|{\rm lcm}(f,n)$ and
$\mathbb Q(\zeta_f) \cap \mathbb Q(\zeta_n,g^{1/n})=\mathbb Q(\zeta_{(f,n)})$
otherwise. In case $b$ is even, we
see that $\Delta\nmid {\rm lcm}(f,n)$ on
using that $n$ is squarefree and we are done. It remains
to deal with the case where $b$ is odd, $n$ is even, $\Delta\nmid n$ and $\Delta|{\rm lcm}(f,n)$. The latter divisibility condition implies
that $\Delta|fn$ and hence $b|n$. Put $\gamma=(-1)^{(b-1)/2}(f,\Delta)$.
Using Lemma \ref{classical} and the fact that $b|n$, it follows that
$\sqrt{\Delta/\gamma}=\sqrt{(-1)^{(b-1)/2}b}\in \mathbb Q(\zeta_n)$.
We now distinguish cases according to
the residue mod $4$ of the odd part of $\Delta$, respectively
the odd part of $(f,\Delta)$.
In each case we check, using Lemma \ref{classical}, that
$\sqrt{\gamma}\in \mathbb Q(\zeta_f)$. It
follows that $\sqrt{\gamma}\not\in \mathbb Q(\zeta_{(f,n)})$, for
otherwise from $\sqrt{\Delta/\gamma},\sqrt{\gamma}\in
\mathbb Q(\zeta_n)$ it would follow that
$\sqrt{\Delta}\in \mathbb Q(\zeta_n)$, contradicting our assumption that
$\Delta\nmid n$.
Using $\sqrt{\Delta/\gamma}\in \mathbb Q(\zeta_n)$,
$\sqrt{\gamma}\in \mathbb Q(\zeta_f)$,
$\sqrt{\gamma}\not\in \mathbb Q(\zeta_{(f,n)})$
and (\ref{graad}), the result is completed. \qed

\vfil\eject
\section{Euler products}
In this section we prove some
results that will help us to write down the Euler
product of the sums encountered
in the proof of the main result (Theorem \ref{main}).
\begin{Prop}
\label{multiplicative}
Let $f,h\ge 1$ be integers. Then the function $w:\mathbb N\rightarrow \mathbb N$
defined by
$$w(k)={k\varphi({\rm lcm}(k,f))\over (k,h)\varphi(f)}$$
is multiplicative. Furthermore,
\begin{itemize}
\item[i)] if $p\nmid h$ and $p\nmid f$, then $w(p)=p(p-1)$
\item[ii)] if $p\nmid h$ and $p|f$, then $w(p)=p$
\item[iii)] if $p|h$ and $p\nmid f$, then  $w(p)=p-1$
\item[iv)] if $p|h$ and $p|f$, then $w(p)=1$
\item[v)] if $h$ is odd, then $w(2)=2$.
\end{itemize}
\end{Prop}
This will help us to prove the following lemma.
\begin{lem}
\label{uiltje}
Let $f,h\ge 1$ be integers
with $1\le a\le f,~(a,f)=1$ and $h$ odd. Let
$\Delta$ be a discriminant of a quadratic number field. Let
$b=\Delta/(f,\Delta)$. Put
$$S(b)=\sum_{n=1,~\Delta|{\rm lcm}(n,f)\atop a\equiv 1({\rm mod~}(f,n))}
^{\infty}{\mu(n)\over w(n)}.$$
Let $S_2(b)$ denote the same sum as $S(b)$ but with the restriction that
$2|n$.
Then $S(b)=0$ if $b$ is even and
$$S(b)={\mu(b)A(a,f,h)\over \prod_{p|b}(w(p)-1)}$$
otherwise. Furthermore, $S_2(b)=-S(b)$.
\end{lem}
{\it Proof}. If $b$ is even, then the summation in $S(b)$ runs over
non squarefree $n$ only and hence $S(b)=0$. Next assume that $b$ is odd.
We have
\begin{eqnarray}
S(b) &=& \sum_{n=1,~b|n/(f,n)\atop a\equiv 1({\rm mod~}(f,n))}^{\infty}
{\mu(n)\over w(n)}=\sum_{d|(a-1,f)}\sum_{n=1,~(f,n)=d\atop b|n/d}
^{\infty}{\mu(n)\over
w(n)}\nonumber\\
&=& \sum_{d|(a-1,f)}{\mu(d)\over w(d)}\sum_{n=1,~(f,n)=1\atop b|n}^{\infty}
{\mu(n)\over w(n)}={\mu(|b|)\over w(|b|)}\sum_{d|(a-1,f)}{\mu(d)\over w(d)}
\sum_{n=1,~(f,nb)=1\atop (b,n)=1}^{\infty}{\mu(n)\over w(n)}.\nonumber
\end{eqnarray}
Now by assumption $\Delta$ is a discriminant and $b$ is odd. This implies
that $(f,b)=1$ and $b$ is squarefree. Thus
\begin{eqnarray}
S(b) &=& {\mu(|b|)\over w(|b|)}\prod_{p|(a-1,f)}(1-{1\over w(p)})
\prod_{p\nmid fb}(1-{1\over w(p)})\nonumber\\
&=& {\mu(|b|)\over w(|b|)}
\prod_{p|(a-1,f)}(1-{1\over w(p)})
\prod_{p\nmid f}(1-{1\over w(p)})\prod_{p|b}(1-{1\over w(p)})^{-1}\nonumber\\
&=& {\mu(|b|)A(a,f,h)\over \prod_{p|b}(w(p)-1)},\nonumber
\end{eqnarray}
where we used that $w(p)>1$ for $p|b$ and
$$\prod_{p|(a-1,f)}(1-{1\over w(p)})\prod_{p\nmid f}(1-{1\over w(p)})=A(a,f,h),$$ an identity immediately obtained on invoking
Proposition \ref{multiplicative}.\\
\indent The latter part of the assertion follows easily on using that
$w(2)=2$. \qed
%\vfill\eject

\section{Proof of the main result}
{\it Proof of Theorem} \ref{main}.\\
Suppose that
$b$ is odd. Note that the discriminant of
$\mathbb Q(\sqrt{\gamma})$ equals $\gamma$.\\

\noindent i) The case
$b$ is odd and $(\gamma/a)=1$. Using Lemma \ref{gemier}, Lemma
\ref{twee}
and the observation above, we find that
$c_a(n)=1$ in case $a\equiv 1({\rm mod~}(f,n))$
and $c_a(n)=0$ otherwise. This together with Theorem \ref{lens} and
Lemma \ref{cycdegree} implies
that
\begin{eqnarray}
\label{star}
\varphi(f)\delta(a,f,g)&=&\sum_{n=1\atop 2\nmid n}^{\infty}{\mu(n)\over w(n)}
+\sum_{n=1,~2|n\atop \Delta\nmid {\rm lcm}(n,f)}^{\infty}{\mu(n)\over w(n)}+
2\sum_{n=1,~2|n\atop \Delta|{\rm lcm}(n,f)}^{\infty}{\mu(n)\over w(n)}.\\
&=&\sum_{n=1}^{\infty}{\mu(n)\over w(n)}+
\sum_{n=1,~2|n\atop \Delta|{\rm lcm}(n,f)}^{\infty}{\mu(n)\over w(n)}
={\rm I}_1+S_2(b)\nonumber,
\end{eqnarray}
where furthermore in each sum we restrict to those
integers $n$ such that
$a\equiv 1({\rm mod~}(f,n))$. Now, on using Proposition \ref{multiplicative},
\begin{eqnarray}
{\rm I}_1&=&\sum_{d|(a-1,f)}\sum_{(f,n)=d}{\mu(n)\over w(n)}=\sum_{d|(a-1,f)}{\mu(d)\over w(d)}\sum_{(f,n)=1}{\mu(n)\over w(n)}\nonumber\\
&=& \prod_{p|(a-1,f)}(1-{1\over w(p)})\prod_{p\nmid f}(1-{1\over w(p)})
=A(a,f,h).\nonumber
\end{eqnarray}
The result follows from
(\ref{star}) on invoking Lemma \ref{uiltje} and Proposition
\ref{multiplicative}.\\

\noindent ii) The case $b$ is odd and $(\gamma/a)=-1$. Using Lemma \ref{gemier},
Lemma \ref{twee}, the observation in the beginning of this proof
and (\ref{star}), we find
$$\varphi(f)\delta(a,f,g)={\rm I}_1+S_2(b)-2\sum_{{2|n,~\Delta\nmid n\atop
\Delta|{\rm lcm}(f,n)}\atop a\equiv 1({\rm mod~}(f,n))}
{\mu(n)\over w(n)}.$$
Now
$$\sum_{2|n,~\Delta|{\rm lcm}(f,n)\atop
a\equiv 1({\rm mod~}(f,n))}{\mu(n)\over w(n)}=
\sum_{{2|n,~\Delta\nmid n\atop \Delta|{\rm lcm}(f,n)}\atop
a\equiv 1({\rm mod~}(f,n))}{\mu(n)\over w(n)}
+\sum_{2|n,~\Delta|n\atop a\equiv 1({\rm mod~}(f,n))}{\mu(n)\over w(n)}.$$
In case $\Delta\equiv 0({\rm mod~}4)$, the latter sum is obviously zero.
In case $\Delta\equiv 1({\rm mod~}4)$, a necessary condition for the latter
sum to be non-zero is that $a\equiv 1({\rm mod~}(f,\Delta))$. This
together with the property of Kronecker sums that
$(b/a_1)=(b/a_2)$ if $a_1\equiv a_2({\rm mod~}|b|)$, shows
that $(\gamma/a)=(\gamma/1)=1$
and hence the latter sum is also zero if
$\Delta\equiv 1({\rm mod~}4)$. It thus follows that
$\varphi(f)\delta(a,f,g)={\rm I}_1+S_2(b)-2S_2(b)={\rm I}_1-S_2(b).$ \\

\noindent iii) The case $b$ is even. As in i) it follows that
$\varphi(f)\delta(a,f,g)={\rm I}_1+S_2(b)$. Now
$I_1=A(a,f,h)$ as we have seen
and $S_2(b)=0$ by Lemma \ref{uiltje}. \qed

%\vfill\eject
\section{Applications}
\label{applic}
Setting $a=f=1$ in Theorem \ref{main} and taking $g\in G\cap
\mathbb Z$ we
obtain Hooley's theorem \cite{Hooley}. Setting $a=1$ we
obtain Theorem 4 of \cite{pie2}. Notice
that of all the progressions mod $f$, the
progression $1({\rm mod~}f)$ is the easiest to deal with, since
trivially $c_1(n)=1$ for every $n$.\\
\indent Lenstra \cite[Theorem 8.3]{L}
gave a sketch of a proof of the following result.
\begin{Thm} {\rm \cite{L}}.
Let $g\in G$, $\Delta$ denote the discriminant of
$\mathbb Q(\sqrt{g})$ and
let $h$ be the largest integer such that $g$ is an $h$-th
power. Then $\delta(a,f,g)=0$ if and only if one of the following holds
\begin{itemize}
\item[i)] $(a-1,f,h)>1$;
\item[ii)] $\Delta|f$ and $({\Delta\over a})=1$ (Kronecker symbol);
\item[iii)] $\Delta|3f$, $3|\Delta$, $3|h$ and $({-\Delta/3\over a})=-1$.
\end{itemize}
\end{Thm}
This result very easily follows from Theorem \ref{main}. We leave it to
the reader to show that if $\delta(a,f,g)=0$, then actually ${\cal P}_{a,f,g}$
is finite.\\
\indent If $S$ is any infinite set of prime numbers, denote by $S(x)$ the number of
primes in $S$
not exceeding $x.$ For given integers $a$ and $f,$ denote by $S(x;f,a)$
the number of primes in $S$ not exceeding $x$ that are congruent to $a$
modulo $f.$ We say that $S$ is weakly uniformly distributed mod $f$
(or WUD mod $f$ for short) if
for every $a$ coprime
to $f,$
$$S(x;f,a)\sim {S(x)\over \varphi(f)},$$
where $\varphi(f)$ denotes Euler's totient.
The progressions $a({\rm mod~}f)$ such that the latter asymptotic equivalence
holds are said to get their fair share of primes from $S.$ Thus
$S$ is weakly uniformly distributed mod $f$ if and only if all the progressions
mod $f$ get their fair share of primes from $S.$
Narkiewicz \cite{N} has written a nice survey on the state of knowledge
regarding the
(weak) uniform distribution of many important arithmetical sequences.
Let $D_g$ denote the
set of natural numbers $f$ such that ${\cal P}_g$ is weakly uniformly
distributed modulo $f$.
\begin{Thm}
\label{equidis}
{\rm \cite{pie2}}. {\rm (GRH)}.
Let $g\in G$ and let
$h$ be the largest integer such that $g$ is an $h$-th
power. Write $g=g_1g_2^2$ with $g_1\in \mathbb Z$ squarefree and
$g_2\in \mathbb Q$.
Assume that not both $g_1=21$ and $(h,21)=7,$ then, assuming GRH,
the set $D_g$ of natural numbers $f$ such that the set of primes $p$
such that $g$ is a primitive root mod $p$ is weakly uniformly distributed
mod $d,$ equals
\begin{itemize}
\item[i)] $\{2^n:n\ge 0\}$ if $g_1\equiv 1({\rm mod~4});$
\item[ii)] $\{1,2,4\}$ if $g_1\equiv 2({\rm mod~4});$
\item[iii)] $\{1,2\}$ if $g_1\equiv 3({\rm mod~4}).$
\end{itemize}
In the
remaining case $g_1=21$ and $(h,21)=7,$
we have $D_g=\{2^n3^m:n,~m\ge 0\}.$
\end{Thm}
Using only a formula for $\delta(1,f,g)$
and Theorem \ref{lens} in some special cases, this result was first deduced
in \cite{pie2}. Using the full force of Theorem
\ref{main}, however, a shorter proof of Theorem
\ref{equidis} can be given.\\

\noindent {\it Proof of Theorem} \ref{equidis}. Put
$S_f:=\{A(a,f,h)~|~1\le a\le f,~(a,f)=1\}.$ Let $b$ and $\gamma$ be
as in Theorem \ref{main}. If $b$ is odd, put $\chi(a)=(\gamma/a)$. Notice
that $\chi$ is a character of $(\mathbb Z/f\mathbb Z)^*$. Notice that the dependence
of $\delta(a,f,g)$ on $a$ comes in only through $\chi(a)$
and $A(a,f,h)$, or rather
the factor $\prod_{p|(a-1,f)}(1-1/p)$ of $A(a,f,h)$.\\
\indent Let us first consider the case where $f=2^m$ for some $m\ge 0$.
Notice that $|S_{2^m}|=1$. If $b$ is even, then $\delta(a,f,g)=A(a,f,h)$
and since $|S_{2^m}|=1$ it follows that $f=2^m\in D_g$. If $b$ is odd,
then it is seen that $\chi$ is the identity if and only if $\Delta$ is odd.
From these two assertions the truth of Theorem \ref{equidis} in case
$f=2^m$ easily follows. Notice that if $f=2^m\not\in D_g$, then the image
of $\delta(.,f,g)$ has cardinality two.\\
\indent It remains to deal with the case where $f$ has
an odd prime factor. First let us
consider the case where $f$ is an odd prime. Then, if $b$ is even,
$$\varphi(f)\delta(1,f,g)=A(1,f,h)\ne A(2,f,h)=\delta(2,f,g)\varphi(f)$$ and we do not have equidistribution. Next
assume that $b$ is odd.
If $f=3$, then a short
calculation shows that a necessary and sufficient
condition for equidistribution
to occur, is that $3|g_1,~(3,h)=1,~\mu(|b|)=-1$, $g_1\equiv 1({\rm mod~}4)$
and that the equation
$\prod_{p|b,~p|h}(p-2)\prod_{p|b,~p\nmid h}(p^2-p-1)=5$ has a solution
with $b$ is odd. Now notice that $g$ is a
solution to this if and only
if $g_1=21$ and $(h,21)=7$. Call such a $g$ exceptional.
Next let $f\ge 5$ (with $f$ a prime). Then
$A(1,f,h)<A(2,f,h)=\cdots=A(f-1,f,h)$. If $\chi$ is the identity, then it
follows from this that $\delta(1,f,g)\ne \delta(a,f,g)$ for every
$1<a<f$ with
$(a,f)=1$. If $\chi$ is not the identity, then it is easily seen that
there exist $1<a_1<a_2<f$, $(a_1a_2,f)=1$, with $\chi(a_1)=-\chi(a_2)$. On
using also that $A(a_1,f,h)=A(a_2,f,h)$ we see that $\delta(a_1,f,g)\ne
\delta(a_2,f,g)$. From this it immediately follows that when
$f$ has an odd prime factor and
$g$ is not exceptional, equidistribution does not occur, for if ${\cal P}_g$ is
WUD mod $f$, then ${\cal P}_g$ must also be WUD mod $\delta$ for every
divisor $\delta$ of $f$. It remains to show that when $g$ is
exceptional we have that ${\cal P}_g$ is WUD
mod $2^m3^n$, for every $m\ge 0$ and $n\ge 1$. This follows
easily from Theorem \ref{main} and the calculation made in the
case $f=3$. \qed\\

\indent Rodier \cite{rod}, in
connection with a coding theoretical problem, conjectured that
the density of the primes in ${\cal P}_2$ such that $p\equiv -1,3$
or $19 ({\rm mod~}28)$ is $A/4$. This would follow if 28 were to be in
$D_2$,  but $28\not\in D_2$ by Theorem \ref{equidis}.
From Theorem \ref{main} it follows
that, under GRH, the density of the set
of primes considered by Rodier is $21A/82$, more
precisely each of the three progressions has density $7A/82$.\\
\indent Let $\cal L$ be the set of odd primes $\ell$ such that there are
infinitely many primes  with
$\ell$ a primitive root mod $p$ and $p$ satisfying $p\equiv \pm 1({\rm mod~}\ell)$.
If $\ell\equiv 1({\rm mod~}4)$ then
$\ell\not\in {\cal L}$ by quadratic reciprocity. Modification of some
of the arguments in \cite{HB} yields that ${\cal L}=\{\ell:\ell\equiv 3({\rm mod~4)}\}$ with at most two primes excluded. That $\cal L$ is non-empty is
used in \cite{sasha} to prove a weaker version of a conjecture of Lubotzky
and Shalev on three-manifolds. By Theorem \ref{main} it follows that,
on GRH, the density of primes $p$ such that
$\ell$ is a primitive root mod $p$ and
$p\equiv \pm 1({\rm mod~}\ell)$, is $(2\ell-1)(\ell-1)A/(\ell^2-\ell-1)$
if $\ell\equiv 3({\rm mod~}4)$ and is zero if $\ell\equiv 1({\rm mod~}4)$.\\
\indent The density of primitive roots in $\mathbb F_p^*$ is $\varphi(p-1)/(p-1)$.
Assuming that the primitive roots are equidistributed over $[1,2,\cdots,p-1]$,
one would perhaps expect that the number $\pi_g(x)$ of primes
$p\le x$ such that $g$ is a primitive root mod $p$, behaves asymptotically
as $\sum_{p\le x}\varphi(p-1)/(p-1)$. This is a well-known and old heuristical
idea.
Using the Siegel-Walfisz theorem for
primes in arithmetic progression, it can be easily shown unconditionally
(see e.g. \cite[Lemma 1]{Stephens}), that the latter
sum is asymptotically equal to
$Ax/\log x$. Comparison with Hooley's theorem \cite{Hooley},
then shows that this heuristic is false
in general (under GRH). Let $h\ge 1$ be as usual the largest integer such
that $g$ is an $h$-th power. A less naive heuristic arises on noting
that $g$ is not a primitive root mod $p$ if $g$ is a square mod $p$ or
$(p-1,h)=1$ and that there are $\varphi(p-1)$ primitive roots amongst
the $(p-1)/2$ non-squares mod $p$. It turns out that this
heuristic is asymptotically
exact (under GRH), even on restricting to primes $p$ in a prescribed
arithmetic progression:
\begin{Thm}
\label{heuro}
{\rm \cite{pie3}. (GRH)}. Let $g\in \mathbb Z\backslash\{-1,0,1\}$
and let $h$ be
the largest integer such that $g$ is an $h$th power. Then
\begin{equation}
\label{laatste}
\pi_g(x;f,a)=2\sum_{{p\le x,~({g\over p})=-1\atop p\equiv a({\rm mod~}f)}
\atop (p-1,h)=1}{\varphi(p-1)\over p-1}+O({x\log \log x\over \log^2x}).
\end{equation}
\end{Thm}
The right hand side in (\ref{laatste}) can be evaluated
unconditionally by an elementary but somewhat lengthy calculation
requiring little beyond the Siegel-Walfisz theorem. The left hand side of
(\ref{laatste}) is of course evaluated in this paper (on GRH). Comparison
of the
main terms in both expressions shows that they are equal and Theorem
\ref{heuro} follows.\\

\noindent Note added July 2003: This paper, written
around 1998, will not be published as its main result will
appear as an application of the main (Galois-theoretical) result of \cite{LMS}. For
a preview of the results of \cite{LMS} the reader is referred to the article
by Stevenhagen \cite{St}.

\medskip\noindent {\footnotesize Faculteit WINS,
Universiteit van Amsterdam, Plantage Muidergracht 24,
1018 TV Amsterdam, The Netherlands.\\
e-mail: {\tt moree@science.uva.nl}\\
homepage: http://web.inter.NL.net/hcc/J.Moree/

\end{document}